\documentclass[a4paper,12pt]{amsart}
\usepackage{amsmath}
\usepackage{amssymb}
\usepackage{mathrsfs}
\usepackage{ifthen}
\usepackage{graphicx}
\usepackage[T1]{fontenc} 

\nonstopmode \numberwithin{equation}{section}
\newtheorem{thm}{Theorem}[section]
\newtheorem{cor}[equation]{Corollary}
\newtheorem{lem}[equation]{Lemma}

\theoremstyle{definition}
\newtheorem{defn}{Definition}[section]

\newtheorem{prob}[equation]{Problem}
\newtheorem{rem}{Remark}[section]

\newcounter{minutes}
\divide\time by 60
\newcounter{hours}
\multiply\time by 60
\addtocounter{minutes}{-\time}

\newcounter {own}
\def\theown {\thesection       .\arabic{own}}

\newenvironment{pf}[1][]{%
	\vskip 3mm
	\noindent
	\ifthenelse{\equal{#1}{}}%
	{{\slshape Proof. }}%
	{{\slshape #1.} }%
}%
{\qed\bigskip}

\newcounter{alphabet}

\newenvironment{Thm}[1][]{\refstepcounter{alphabet}%
	\bigskip%
	\noindent%
	{\bf Theorem \Alph{alphabet}}%
	\ifthenelse{\equal{#1}{}}{}{ (#1)}%
	{\bf .} \itshape}{\vskip 8pt}

\def\be{\begin{equation}}
	\def\ee{\end{equation}}

\newcommand{\bee}{\begin{enumerate}}
	\newcommand{\eee}{\end{enumerate}}

\newcommand{\blem}{\begin{lem}}
	\newcommand{\elem}{\end{lem}}
\newcommand{\bthm}{\begin{thm}}
	\newcommand{\ethm}{\end{thm}}
\newcommand{\bcor}{\begin{cor}}
	\newcommand{\ecor}{\end{cor}}
\newcommand{\beg}{\begin{examp}}
	\newcommand{\eeg}{\end{examp}}
\newcommand{\begs}{\begin{examples}}
	\newcommand{\eegs}{\end{examples}}
\newcommand{\bdefe}{\begin{defin}}
	\newcommand{\edefe}{\end{defin}}
\newcommand{\bprob}{\begin{prob}}
	\newcommand{\eprob}{\end{prob}}
\newcommand{\bei}{\begin{itemize}}
	\newcommand{\eei}{\end{itemize}}

\begin{document}
	
	\title{Neargeodesics in Gromov hyperbolic John domains in Banach spaces}

	\author{Vasudevarao Allu}
	\address{Vasudevarao Allu,
		Discipline of Mathematics,
		School of Basic Sciences,
		Indian Institute of Technology  Bhubaneswar,
		Argul, Bhubaneswar, PIN-752050, Odisha (State),  India.}
	\email{avrao@iitbbs.ac.in}

	\author{Abhishek Pandey}
	\address{Abhishek Pandey,
		Discipline of Mathematics,
		School of Basic Sciences,
		Indian Institute of Technology  Bhubaneswar,
		Argul, Bhubaneswar, PIN-752050, Odisha (State),  India.}
	\email{ap57@iitbbs.ac.in}

	\subjclass[2010]{Primary 30C65, 30L10, 30F45. Secondary 30C20}
	\keywords{Quasihyperbolic metric, Gromov hyperbolic spaces, John domain, neargeodesics, quasiconformal mappings, quasihyperbolic geodesic, cone arc}

	\def\thefootnote{}
	\footnotetext{ {\tiny File:~\jobname.tex,
			printed: \number\year-\number\month-\number\day,
			\thehours.\ifnum\theminutes<10{0}\fi\theminutes }
	} \makeatletter\def\thefootnote{\@arabic\c@footnote}\makeatother

	\thanks{}
	
	\maketitle
	\pagestyle{myheadings}
	\markboth{Vasudevarao Allu and Abhishek Pandey}{Neargeodesics in Gromov hyperbolic John domains in dimensional Banach spaces}
	
	\begin{abstract}
	In this paper, we prove that neargeodesics in Gromov hyperbolic John domains in Banach space are cone arcs. This result gives an improvement of a result of Li [Theorem 1, Int. J. Math. 25 (2014)].
	\end{abstract}

\section{Introduction and main result}	
Let $\Omega$ be an open ball or half space in $\mathbb{R}^n$, we denote $h_{\Omega}$ be the hyperbolic metric in $\Omega$ with constant curvature $-1$. If $\Omega$ is an upper half plane then $h_{\Omega}$ is defined in $\Omega$ by the means of the density
$$\rho(z)=\frac{1}{{\rm Im \,z}}=\frac{1}{d(z,\partial \Omega)}.$$
This density can be used to introduce an analogue of the hyperbolic metric in an arbitrary subdomain $D$ in $\mathbb{R}^n$. The quasihyperbolic metric in $\mathbb{R}^n$  is a generalization of hyperbolic metric, which was introduced by Gehring and Palka \cite{Gehring-1976} and they have studied the characterization of domains $D$ in the one point compactification of $\mathbb{R}^n$ which are quasiconformally equivalent to ball in $\mathbb{R}^n$. Infact, such domains are homogeneous via quasiconformal mappings. Gehring and Palka \cite{Gehring-1976} have proved that the maximal dilatation of such quasiconformal mappings can be estimated in terms of quasihyperbolic metric and using these estimates they obtained useful results related to the homogeneity of domains with respect to quasiconformal family. Thus, the importance of the quasihyperbolic metric is quite clear as it is well-behaved with the quasiconformal mappings. The quasihyperbolic metric can be naturally defined on a wide range of metric spaces, including Banach spaces. In this case quasiconformal maps are defined in terms of quasihyperbolic metric and the study of this concept is known as free quasiconformality which has been introduced and developed by V\"ais\"al\"a (see \cite{Vaisala-1990,Vaisala-1991,Vaisala-1999}). 
\par In this paper, our main focus is to study the geometric properties of quasihyperbolic quasigeodesic in John domains in Banach space.\\

It is well known that, in the Poincar\'e disk model of the hyperbolic space, hyperbolic geodesic between two points $x,y\in \mathbb{D}$, denote it by $\gamma_{hyp}[x,y]$, have the following properties:
\begin{itemize}
	\item[(i)] $l(\gamma_{hyp}[x,y])\le C|x-y|$, and
	\item[(ii)] $\min\{l(\gamma_{hyp}[x,z]),l(\gamma_{hyp}[z,y])\}\le C\,d(z,\partial D)$,
\end{itemize}
for all $z\in \gamma$, where $C=\pi/2$. The second condinition is known as the cone arc condition by which we defined the John domains, see Definition \ref{John}, which was first introduced by John \cite{John-1961} in the context of elasticity theory. If the curve satisfies both conditions mentioned above it is called double cone arc by which we define the uniform domains, see Definition \ref{Uniform}. Gehring and Osgood \cite{Gehring-1979} have proved that each quasihyperbolic geodesic in a $c$-uniform domain $D\in \mathbb{R}^n$ is a double $b$-cone arc, where the constant $b$ depends only on $c$. Since John domains can be thought of one sided uniform domains, it is natural to ask whether such result also holds for John domains or not? In fact, this problem has been proposed by Gehring {\it et. al.} \cite{Gehring-1989} in the following form.
\begin{prob}\label{Problem}
	Suppose $D\subset \mathbb{R}^n$ is a $c-$ John domain and that $\gamma$ is a quasihyperbolic geodesic in $D$. Is $\gamma$ a $b$-cone arc for some $b=b(c)$?
\end{prob}

In 1989, Gehring {\it et. al.} \cite{Gehring-1989} proved the following result.

\begin{Thm}\cite[Theorem 4.1]{Gehring-1989}
If $D\subset \mathbb{R}^2$ is a simply connected $c$-John domain then every quasihyperbolic or hyperbolic geodesic in $D$ is a $b$-cone arc, where $b$ only depends on $c$.
\end{Thm}

In fact,  Gehring {\it et. al.} \cite{Gehring-1989} have constructed several examples that shows that a quasihyperbolic geodesic in a $c$-John domain need not be a $b$-cone arc with $b=b(c)$ unless $n=2$ and $D$ is simply connected. Thus, this suggests us that we need some extra condition on the $c$-John domain so that the answer to the Problem \ref{Problem} is affirmative. 
\par This raises the following natural problem:
\begin{prob}\label{Big}
Determine necessary and/or sufficient conditions for quasihyperbolic geodesics in a John space to be cone arcs.
\end{prob}

In $1989$, Heinonen \cite[Question 2]{Heinonen-1989} posed the following problem.
\begin{prob}\label{Problem-1}
	Suppose $D\subset \mathbb{R}^n$ is a bounded $c$-John domain, that is quasiconformally equivalent to the unit ball $\mathbb{B}^n$ in $\mathbb{R}^n$, and $\gamma$ is a quasihyperbolic geodesic in $D$. Is $\gamma$ a $b=b(c)$-cone arc for some constant $b$? 
\end{prob}

In $2001$, Bonk, Heinonen, and Koskela \cite[Proposition 7.12]{Koskela-2001} gave an affirmative answer to Probelm \ref{Problem-1}.  To be more specific, every quasihyperbolic geodesic in a bounded Gromov hyperbolic John domain in $\mathbb{R}^n$ is a cone arc, as can be seen in (\cite[Theorem B]{Li-2014} or \cite[Page 228]{Rasila-2022}). Note that their cone arc constant depends on the space dimension $n$.
\par Recently, various contributions have been made in this line. In 2022, Zhou, Li, and Rasila \cite{Rasila-2022} consider the problem \ref{Big} in the setting of metric space and obtained a dimension-free answer to the Problem \ref{Problem-1}. Precisely, they proved the following: 

\begin{Thm}\cite[Theorem 1.1]{Rasila-2022}\label{6-1}
	Let $(D,d)$ be a locally compact, rectifiably connected and non complete metric space, and $k$ be the quasihyperbolic metric on $(D,d)$. Suppose $(D,d)$ is $a$-John. If $(D,k)$ is $K$-roughly starlike and Gromov $\delta$-hyperbolic, then every quasihyperbolic geodesic in $D$ is $b$-cone arc, where $b$ depends only on $a,K$ and $\delta$.
\end{Thm}

\begin{rem}\cite[Corollary 1.2]{Rasila-2022}
	Observe that any proper subdomain $D$ in $\mathbb{R}^n$ is locally compact, rectifiably connected and noncomplete with respect to the usual metric. Moreover, every $\delta$-hyperbolic domain is roughly starlike for some $K$ depending only on $\delta$ (see \cite[Theorem 3.22]{Vaisala-GH-1}).
\end{rem}

The next result of Zhou and Ponnusamy \cite{Ponnusamy-2024-I} shows that the condition of rough starlikeness is not required in Theorem B. Thus, the following result of Zhou and Ponnusamy \cite{Ponnusamy-2024-I} gives an improvement of Theorem B.

\begin{Thm}\cite[Theorem 1.3]{Ponnusamy-2024-I}\label{6-2}
	Let $(D,d)$ be a locally compact, rectifiably connected and non complete metric space, and $k$ be the quasihyperbolic metric on $(D,d)$. Suppose $(D,d)$ is $a$-John space. If $(D,k)$ is Gromov $\delta$-hyperbolic, then every quasihyperbolic geodesic in $D$ is a $b$-cone arc with $b=b(a,\delta)$.
\end{Thm}
The main aim of this paper is to study Problem \ref{Big} in the Banach space setting. Note that the quasihyperbolic geodesic may not exists in the infinite-dimensional Banach spaces (see \cite[Example 2.9]{Vaisala-1990} and \cite[Remark 3.5]{Vaisala-1999}). To overcome this shortage, V\"ais\"al\"a \cite{Vaisala-1991} introduced the concept of neargeodesics. We recall the definition of neargeodesic.
\begin{defn}
Let $E$	be a real Banach space and $D\subsetneq E$ be a domain. Let $c\ge 1$. A curve $\gamma$ in $D$ is said to be $c$-neargeodesic if for all $x,y \in \gamma$, we have 
	$$l_k(\gamma[x,y])\le c\,k_D(x,y).$$
\end{defn}
Thus a curve $\gamma$ is a quasihyperbolic geodesic if, and only if, it is a $1$-neargeodesic. Neargeodesics can also be termed as quasihyperbolic quasigeodesics. In \cite{Vaisala-1991} V\"ais\"al\"a proved the following existence theorem for neargeodesics. 

\begin{Thm}\cite[Theorem 3.3]{Vaisala-1991}
	Let $c>1$. Then for every pair of points $z_1,z_2\in D$ there exists a $c$-neargeodesic joining $z_1$ and $z_2$.
\end{Thm}

In 2014, Y. Li \cite[Theorem 1]{Li-2014} considered the Problem \ref{Problem-1} further in the Banach space setting and proved the following

\begin{Thm}\cite[Theorem 1]{Li-2014}\label{Li}
	Suppose that $D\subset E$ is an $a$-John domain which is homeomorphic to an $c$-inner uniform domain via an $(M, C)$-CQH. Let $z_1,z_2\in D$ and $\gamma$ is a $c_0$-neargeodesic joining $z_1$ and $z_2$ in $D$, then $\gamma$ is a $b$-cone arc with $b=b(a,c,c_0, M, C)$.
\end{Thm}
\begin{rem}
\begin{itemize}
	\item[(1)] A homeomorphism $f:D\to D'$ is said to be $C$-coarsely, $M$-quasihyperbolic, or briefly $(M,C)$-CQH if it satisfies
	$$\frac{k_D(x,y)}{M}-C\le k_{D'}(f(x),f(y))\le M k_{D}(x,y)+C$$
	for all $x,y\in D$.
	\item[(2)] Observe that the John domain in Theorem E is Gromov hyperbolic since it is an image of a Gromov hyperbolic domain (an inner uniform domain) under a quasihyperbolic $(M,C)$  quasi-isometry). 
\end{itemize}
\end{rem}

This motivates us to ask does Theorem C type result still hold if we consider John domains in Banach space? If the answer to this question is yes, then it will give an improvement of Theorem E of Y. Li. Indeed, using the ideas of \cite{Ponnusamy-2024-I}, we prove that the answer is yes. Our main result is as follows.

\begin{thm}\label{main}
	Let $E$ be a real Banach space and $D\subsetneq E$ be a $c$-John domain. If $D$ is Gromov hyperbolic, {\it i.e.} $(D,k_D)$ is $\delta$-hyperbolic for some $\delta\ge 0$, then every $c_0$-neargeodesic in $D$ is $b$-cone arc, where $b$ depends only on $c$, $c_0$ and $\delta$.
\end{thm}

	\section{Preliminaries}
	In this section we mention some basic as well as advanced concepts and related results in the literature, which are useful to prove our main result.
	
	\subsection{Metric Geometry}
	Let $(X,d)$ be a metric space. A curve is a continuous function $\gamma:[a,b]\rightarrow X$. If a curve $\gamma:[a,b]\rightarrow X$ is an  embedding of $[a,b]$, then it is called an arc.
	Let $\mathcal{P}$ denote set of all partitions $a = t_{0}<t_{1}<t_{2}< \cdots<t_{n}=b$ of the interval $[a,b]$. The  length of the curve $\gamma$ in the metric space $(X,d)$ is 
	$$l_{d}(\gamma) = \sup_{\mathcal{P}} \sum_{k=0}^{n-1}d(\gamma(t_{k}), \gamma(t_{k+1})).$$
	A curve is said to be rectifiable if $l_{d}(\gamma) < \infty$.  A metric space $X$ is said to be rectifiably connected if every pair of points $x,y \in X$ can be joined by a rectifiable curve. For a rectifiable curve $\gamma$ we define arc length $s:[a,b]\rightarrow \left[ 0,l_d(\gamma) \right] $ by $s(t) = l_{d}(\gamma|_{[a,t]})$. The arc length function is of bounded variation. For any rectifiable curve $\gamma:[a,b]\to X$, there is a unique map $\gamma_s:[0,l_d(\gamma)]\to X$ such that $\gamma=\gamma_s\circ s$, and such a curve $\gamma_s$ is called the arclength parametrization of $\gamma$.\\
	\par For $x,y\in X$, the inner length metric $\lambda_X(x,y)$ is defined by
	$$\lambda_X(x,y)=\inf\{l_d(\gamma): \gamma \mbox{ is a rectifiable curve joining } x \mbox{ and } y\}.$$
	
	Let $\rho :X\rightarrow[0,\infty]$ be a Borel function. The $\rho$-length of a rectifiable curve $\gamma$ is $$\int_{\gamma}\rho\,ds=\int_{a}^{b}\rho(\gamma(t))\, ds(t)=\int_{0}^{l_d(\gamma)}\rho\circ \gamma_s(t)\,dt.$$
	
	If $X$ is rectifiably connected then $\rho$ induces a distance function which is defined by $$d_{\rho}(x,y) = \inf \int_{\gamma}\rho\,ds,$$
	where infimum is taken over all rectifiable curves joining $x$ and $y$ in $X$. We note that, in general, $d_\rho$ need not to be a metric however $d_{\rho}$ is a metric if $\rho$ is positive and continuous. If $d_\rho$ is a metric, we say $(X,d_\rho)$ is a conformal deformation of $(X,d)$ by the conformal factor $\rho$.

		\begin{defn}
		A curve $\gamma:[a,b]\to X$ is said to be geodesic if for all $t,t'\in [a,b]$, $$d(\gamma(t), \gamma(t'))=|t-t'|.$$
	 A metric space $(X,d)$ is said to be
		\begin{itemize}
			\item[(1)] {\it a geodesic space} if every pair of points $x,y \in X$ can be joined by a geodesic,
			\item[(2)] {\it proper} if every closed ball is compact in $X$,
			\item[(3)] {\it intrinsic} or a {\it path metric space} or a {\it length metric space} if for all $x,y\in X$, it holds that
			$$d(x,y)=\lambda_X(x,y).$$ 	
		\end{itemize}
	\end{defn}
	
	\subsection{Gromov hyperbolicity} In the 1980s, Gromov established the concept of Gromov hyperbolicity. The underlying space is frequently taken to be a geodesic metric space, and in most cases, it is a proper metric space. We adapt the concept of Gromov hyperbolicity, which is effective for non-geodesic spaces, since, as we shall see later, an infinite-dimensional Banach space with quasihyperbolic metric need not be a geodesic space. For this, we refer to the work of V\"ais\"al\"a \cite{Vaisala-GH}, where the Gromov product concept is applied.
	
	\begin{defn}
		Let $(X,d)$ be a metric space and let $p \in X$. The Gromov Product $(x|y)_p$ of $x,y \in X$ with respect to $p$ is defined by 
		$$ (x|y)_{p} = \frac{1}{2}\Big(d(x,p)+d(y,p)-d(x,y)\Big).$$
		Therefore, $(x|y)_{p}$ measures to what extent the triangle inequality for the triangle $\Delta xyp$ differs from equality.
	\end{defn}

	\begin{defn}
		Let $\delta\ge 0$. A metric space $(X,d)$ is said to be $\delta$-hyperbolic if for every $x,y,z,p \in X$, the following inequality holds
		$$(x|z)_{p} \geq \min \Bigl\{(x|y)_{p},(y|z)_{p}\Bigr\} - \delta.$$
		A metric space $(X,d)$ is said to be Gromov hyperbolic if it is $\delta$-hyperbolic for some $\delta\ge 0$.
	\end{defn}
	
	\subsection{Examples} Clearly, $\mathbb{R}$ is a 0-hyperbolic. Any tree $T$ is $0$-hyperbolic. The unit disk with hyperbolic metric is Gromov hyperbolic with $\delta=\log\,3$. The complex plane is not Gromov hyperbolic.\\
	
	We shall mention an important result here due to V\"ais\"al\"a \cite[3.7]{Vaisala-GH}, which is known as the stability theorem. In fact, the stability theorem says that in an intrinsic hyperbolic space any two quasi-isometric paths with same initial and terminal points $x$, $y$ run close to each other even if $|x-y|$ is large. To mention this theorem, we first need to mention the concepts of quasi-isometry and Hausdorff distance.
	\begin{defn}
		Let $X$ and $Y$ be two non-empty subsets of a metric space $(M,d)$. We define their Hausdorff distance $d^{\mathcal{H}}(X,Y)$ by 
		$$d^{\mathcal{H}}(X,Y)=\max\left\{\sup_{x\in X}d(x,Y),\, \sup_{y\in Y}d(X,y)\right\},$$
		where $d(a,B)=\inf_{b\in B}d(a,b)$ denotes the distance from a point $a\in M$ to a subset $B\subset M$.
	\end{defn}
	
	\begin{defn}
		Let $\lambda\ge 1$ and $\mu\ge 0$ and $(X,d)$ and $(Y,d')$ be metric spaces. We say that a map $f:X\to Y$ is a $(\lambda,\mu)$-quasi-isometry if
		$$\lambda^{-1}d(x,y)-\mu\le d'(f(x),f(y))\le \lambda d(x,y)+\mu$$
		for all $x,y\in X$. In the case where $f:I\to Y$ is a map of a real interval $I$, we say that such a map $(\lambda,\mu)$-quasi-isometric path.
	\end{defn}
	
	\begin{rem}\label{rem-1}
		An arc $\gamma$ joining $x$ and $y$ in a metric space $X$ is said to be $\lambda$-quasigeodesic, $\lambda\ge 1$, if 
		$$l(\gamma[u,v])\le \lambda d(u,v)$$
		for all $u,v\in \gamma$. In such a case the arc length parametrization $\gamma_s:[0,l(\gamma)]\to \gamma$ satisfies the inequality
		$$\lambda^{-1}|t-t'|\le d(\gamma_s(t),\gamma_s(t'))\le \lambda |t-t'|.$$
		Thus $\gamma$ is $(\lambda,0)$-quasi-isometry.
	\end{rem}

	The following theorem (see \cite[3.7]{Vaisala-GH}) is vital to prove our main result.
	
	\begin{Thm}\label{HD}
		Suppose that $X$ is an intrinsic $\delta$-hyperbolic space and that $\gamma$ and $\alpha$ are $(\lambda,\mu)$-quasi-isometric path with same initial and terminal points. Then there exists a constant $M$ such that $$d^{\mathcal{H}}(\gamma, \alpha)\le M,$$ where $M$ depends only on $\delta,\lambda,\mu$ and $d^{\mathcal{H}}$ denote the Hausdorff distance between $\gamma$ and $\alpha$.
	\end{Thm}
	
	Next we mention an important result due to V\"ais\"al\"a (see \cite[Theorem 3.18]{Vaisala-GH}) which says that quasi-isometries preserves hyperbolicity. 
	
	\begin{Thm}\label{hyperbolic}
		Let $X$ and $Y$ be intrinsic spaces and let $f: X \to Y$ be a $(\lambda,\mu)$-quasi-isometry. If $Y$ is $\delta$-hyperbolic, then $X$ is $\delta'$-hyperbolic with $\delta'=\delta'(\delta,\lambda,\mu)$.
	\end{Thm}

	\subsection{Quasihyperbolic metric} In view of our frame work, we consider the quasihyperbolic metric in the setting of Banach spaces. 
	\begin{defn}\label{Def-1}
		Let $E$ be a real Banach space of dimension at least $2$ and let $D\subsetneq E$ be a domain (open, connected nonempty set). Let the metric induced by norm on D is $d$. Define $k: D \to (0,\infty)$ by
		$$k(z)=\frac{1}{\Delta(z)}=\frac{1}{d(z,\partial D)}.$$
		Define the quasihyperbolic length of a rectifiable curve $\gamma$ in $D$ by
		$$l_{k}(\gamma)=\int_{\gamma}\rho(z)\,ds=\int_{\gamma}\frac{ds}{d(z,\partial D)}.$$
		Let $k_D$ denote the quasihyperbolic metric in $D$ which is defined by
		\begin{equation}\label{qh metric}
			k_{D}(x,y)=\inf_{\gamma}l_{k}(\gamma)
		\end{equation}
		where the infimum is taken over all rectifiable curves $\gamma$ in $D$ joining $x$ to $y$.
	\end{defn}
	Observe that $(D,k_D)$ is always an intrinsic space. For further geometric properties of the quasihyperbolic metric, we refer to \cite{Gehring-1979} and  \cite{Matti-book}. Quasihyperbolic metric satisfies the following properties.
	
	\begin{itemize}
		\item[(i)] $k_{D}=h_{D}$ when $D$ is a half space in $\mathbb{R}^n$.
		\item[(ii)] $k_{D}\le h_{D}\le 2k_{D}$ when $D$ is a ball in $\mathbb{R}^n$.
	\end{itemize}

Finally, we recall some estimates for the quasihyperbolic metric which have been first introduced by Gehring and Palka \cite[2.1]{Gehring-1976} in $\mathbb{R}^n$. Later, V\"ais\"al\"a \cite[Lemma 2.2]{Vaisala-1990} have proved these inequalities in the case of Banach spaces. Let $D\subsetneq E$ be a domain and $x,y\in D$ and let $\gamma$ be a rectifiable curve joining $x$ and $y$. Then we have the following:
	\begin{equation}\label{qh-eq-1}
		k(x,y)\ge \log\left(1+\frac{|x-y|}{\min\{\Delta(x),\Delta(y)\}}\right)\ge \left|\log \frac{\Delta(y)}{\Delta(x)}\right|
	\end{equation}
	and 
	\begin{equation}\label{qh-eq-2}
		l_k(\gamma)\ge \log\left(1+\frac{l(\gamma)}{\min\{\Delta(x),\Delta(y)\}}\right).
	\end{equation}

	We are going to frequently use  \eqref{qh-eq-1} and \eqref{qh-eq-2} to prove our main results. 
	
	\subsection{Quasihyperbolic geodesic}
	A rectifiable curve $\gamma$ from $a$ to $b$ in $D$ is said to be quasihyperbolic geodesic if $k_{D}(a,b)=l_{k}(\gamma)$.
	Obviously each subarc of a quasihyperbolic geodesic is again a geodesic. Observe that quasihyperbolic geodesic is curve $\gamma$ for which the infimum in (\ref{qh metric}) is attained.\\
	
At this juncture, the natural question is: does quasihyperbolic geodesic always exist? A result of Gehring and Osgood \cite{Gehring-1979} shows that for a proper subdomain $D$ in $\mathbb{R}^n$, the answer to this question is affirmative. For important results on quasihyperbolic geodesics in domain in $\mathbb{R}^n$, we refer to G. Martin \cite{Martin-1985}. Moreover, Martin \cite{Martin-1985} has proved that quasihyperbolic geodesics in $\mathbb{R}^n$ are $C^{1,1}$ smooth, {\it i.e.} the arc length parametrization has Lipschitz continuous derivatives. It is well known that a quasihyperbolic geodesic between any two points exists if dim ($E$) is finite (see \cite[Lemma 1]{Gehring-1979}). Note that the quasihyperbolic geodesic may not exist in infinite-dimensional Banach spaces (see \cite[Example 2.9]{Vaisala-1990} and \cite[3.5]{Vaisala-1999}). These examples of domains are in Hilbert spaces, such that the complements of these domains are uncountable. For an example of a domain that is not geodesic with respect to quasihyperbolic metric such that the complement of the domain is countable, we refer to \cite[Example 4.1]{Rasila-2017}. However, V\"ais\"al\"a \cite[Theorem 2.1]{Vaisala-2005} has proved the existence of quasihyperbolic geodesic when $E$ is reflexive Banach space and $D$ is a convex domain. Furthermore, in 2011, Martio and V\"ais\"al\"a \cite[Theorem 2.11]{Martio-2011} proved the existence as well as uniqueness of a quasihyperbolic geodesic when $E$ is a uniformly convex Banach space and $D$ is a convex domain. 
	
	
	\subsection{Uniform Domains}
	
	\begin{defn}\label{Uniform}
		Let $c\ge 1$. A domain $D$ in $E$ is said to be $c$-uniform in the norm metric if for each pair of points $x,y\in D$, there exists a rectifiable arc $\gamma$ in $D$ joining $x$ to $y$ such that 
		\begin{itemize}
			\item[(i)] $\min\{l(\gamma[x,z]),l(\gamma[z,y])\}\le c\,d(z,\partial D)$, for all $z\in \gamma$, and
			\item[(ii)]	$l(\gamma)\le c|x-y|$,	
		\end{itemize}
		where $l(\gamma)$ denotes the arc length of $\gamma$ in $(D,|.\,|)$, where $|.\,|$ is metric induced by norm, and $\gamma[x,z]$ denotes the part of $\gamma$ between $x$ and $z$. Also we say that curve $\gamma$ is a double $c$-cone arc.
		
	\end{defn}
	
	\subsection{Inner uniform domains}
	\begin{defn}
		Let $c\ge 1$. A domain $D$ is said to be $c$-inner uniform in the norm metric if for each pair of points $x,y\in D$, there exists a rectifiable cuve $\gamma$ in $D$ joining $x$ to $y$ such that 
		\begin{itemize}
			\item[(i)] $\min\{l(\gamma[x,z]),l(\gamma[z,y])\}\le c\,d(z,\partial D)$, for all $z\in \gamma$,
			\item[(ii)]	$l(\gamma)\le c \lambda_D(x,y)$.
		\end{itemize}	
		
	\end{defn}
	
	\subsection{John Domain}

	\begin{defn}\label{John}
		Let $c\ge 1$. A domain $D$ is said to be $c$-John in the norm metric if for each pair of points $x,y\in D$ there exists a rectifiable curve $\gamma$ in $D$ joining $x$ to $y$ such that for all $z\in \gamma$ 
		\begin{equation}\label{John-2}
			\min\{l(\gamma[x,z]),l(\gamma[z,y])\}\le c\,d(z,\partial D)
		\end{equation}
		and we say that such a curve $\gamma$ is a $c$-cone arc.
	\end{defn}
	
	\begin{rem} We mention some of the remarks with regards to the connection between the aforesaid domains.
		
		\begin{itemize}
			\item[(i)] $c$-uniform implies inner $c$-uniform implies $c$-John
			\item[(ii)] $\mathbb{R}^2\setminus \{(x,0):x\ge 0\}$ is inner uniform but not uniform.
			\item[(iii)] $\mathbb{R}^2\setminus \{(n,0):n\in \mathbb{N}\}$ is a John domain but not inner uniform.
		\end{itemize}
	\end{rem}
	
	\subsection{Gromov hyperbolicity of domains}
	\begin{defn}
		Let $\delta\ge 0$. A domain $D$ in a Banach space $E$ is said to be $\delta$-hyperbolic if $(D,k_D)$ is $\delta$-hyperbolic. 
	\end{defn}
	
	\begin{defn}
		A domain $D$ in a Banach space $E$ is said to be Gromov hyperbolic if $(D,k_D)$ is $\delta$-hyperbolic for some $\delta\ge 0$.
	\end{defn}
	
	All $c$-uniform and $c$-inner uniform domains are Gromov $\delta=\delta(c)$ hyperbolic. $\mathbb{R}^2\setminus \{(n,0):n\in \mathbb{N}\}$ is a John domain which is not hyperbolic. The broken tube (see \cite[2.12]{Vaisala-1990} and for a detailed treatment we refer to \cite{Vaisala-2004}) in an infinite dimensional seperable Hilbert space is a classical example of a domain which is hyperbolic but not John and hence not uniform. 
	
	\section{Proof of Theorem \ref{main}}
	We shall divide the proof of Theorem \ref{main} into two theorems. Theorem \ref{thm-1} below gives a sufficient condition for a neargeodesic to be a cone arc and Theorem \ref{thm-2} below says that in Gromov hyperbolic John domains, every neargeodesic satisfies this sufficient condition. Therefore, combining these two theorem we obtain the proof of our main result.
	
	\begin{thm}\label{thm-1}
		Let $\gamma$ be a $c_0$- neargeodesic joining $x,y$ such that the following holds: If $z_1,z_2 \in \gamma$ such that each $w\in \gamma[z_1,z_2]$ satisfies $\Delta(w)\le 2\min\{\Delta(z_1), \Delta(z_2)\}$, then there exists a constant $A\ge 1$ such that  $$k_D(z_1,z_2)\le A.$$ 
		Then $\gamma$ is a $C$-cone arc, where $C$ depends only on $A$ and $c_o$.
	\end{thm}
	
	\begin{pf}
		Suppose $\gamma$ is a $c_0$- neargeodesic joining $x,y$ and if $z_1$ and $z_2$ are points on $\gamma$ satisfying for each $w\in \gamma[z_1,z_2]$, $\Delta(w)\le 2\min\{\Delta(z_1), \Delta(z_2)\}$, then there exists a constant $A\ge 1$ such that  $$k_D(z_1,z_2)\le A.$$ 
		Since $\Delta$ is continous and $\gamma$ is compact, choose a point $z_0\in \gamma$ with 
		$$\Delta(z_0)=\max_{z\in \gamma}\Delta(z)$$
		Then there exists a unique non-negative integer $n$ such that
		$$2^n\Delta(x)\le \Delta(z_0)\le 2^{n+1}\Delta(x).$$
		Let $x=x_0$ and for each $i=1,2,\ldots,n$, let $x_i$ be the first point on $\gamma[x,z_0]$ such that
		$$\Delta(x_i)=2^i \Delta(x_0).$$
		Similarly, there exists $m\ge 0$ such that
		$$2^m\Delta(y)\le \Delta(z_0)\le 2^{m+1}\Delta(y).$$
		Let $y=y_0$ and for each $j=1,2,\ldots,m$, let $y_j$ be the last point on $\gamma[z_0,y]$ such that
		$$\Delta(y_j)=2^j \Delta(y_0).$$
		In this manner we have devided the arc $\gamma$ into $(n+m+1)$ subarcs 
		$$\gamma[x_0,x_1],\ldots, \gamma[x_{n-1},x_n], \gamma[x_n,y_m],\gamma[y_m,y_{m+1}],\ldots,\gamma[y_1,y_0].$$
		It is easy to observe that all these subarcs are also $c_0$-neargeodesic between their respective end points. Let $u$ and $v$ be any two adjacent points along $\gamma$. That is $u$ and $v$ are $x_i, x_{i+1}$ or $x_n, y_m$ or $y_{j+1}, y_j$, $0\le i\le n-1$ and $0\le j\le m-1$. By the construction, it is clear that for each $z\in \gamma[u,v]$, we have
		$$\Delta(z)\le 2\min\{\Delta(u), \Delta(v)\}.$$
		By assumption, there exists a constant $A\ge 1$ such that 
		$$k_D(u,v)\le A.$$
		To show that $\gamma$ is a $C$-cone arc for some $C$, we need to show that
		\begin{equation}\label{cone}
			\min\{l(\gamma[x,w]), l(\gamma[w,y])\}\le C\Delta(w)\,\,\,\,\, \mbox{ for all } w\in \gamma.
		\end{equation}
		In order to prove \eqref{cone}, it is sufficient to prove either $l(\gamma[x,w])\le C\Delta(w)$ or $l(\gamma[w,y])\le C\Delta(w)$. For this we need to consider three cases depending upon the location of $w$ on $\gamma$. 
		\begin{itemize}
			\item[\textbf{Case 1.}] Suppose $w\in \gamma[x_i,x_{i+1}]$ for some $0\le i\le n-1$. Then in view of \eqref{qh-eq-1}, we have
			$$\log\frac{\Delta(x_i)}{\Delta(w)}\le k_D(w,x_i)\le l_k(\gamma[x_i,w])\le l_k(\gamma[x_i,x_{i+1}]) \le c_0\,k_D(x_i,x_{i+1})\le c_0\,A.$$
			Therefore, we have 
			\begin{equation}\label{equ-1}
				\Delta(x_i)\le e^{c_0A}\Delta(w).
			\end{equation}
			Moreover, we have
			\begin{equation}\label{equ-2}
				l(\gamma[x,w])\le \sum_{s=0}^{i}l(\gamma[x_s,x_{s+1}])
			\end{equation}

			In view of \eqref{qh-eq-2}, we have
			\begin{eqnarray*}
				\log\left(1+\frac{l(\gamma[x_s,x_{s+1}])}{\min\{  	\Delta(x_s),\Delta(x_{s+1})\}}\right)&\le& l_k(\gamma[x_s,x_{s+1}])\\ \nonumber
				&\le& c_0\,k_D(x_s,x_{s+1})\\ \nonumber
				&\le& c_0\,A
			\end{eqnarray*}
			Therefore, this gives us that
			$$\log\frac{l(\gamma[x_s,x_{s+1}])}{\Delta(x_s)}\le \log\left(1+\frac{l(\gamma[x_s,x_{s+1}])}{\Delta(x_s)}\right)\le c_0\,A$$
			and hence
			\begin{equation}\label{equ-3}
				l(\gamma[x_s,x_{s+1}])\le e^{c_0\,A}\Delta(x_s)
			\end{equation}
			Using \eqref{equ-3} in \eqref{equ-2}, we obtain
			\begin{eqnarray}
				l(\gamma[x,w])&\le& \sum_{s=0}^{i}e^{c_0\,A}\Delta(x_s)\\ \nonumber
				&\le&2\,e^{c_0\,A}\,\Delta(x_i)\\ \nonumber
				&\le&2\,e^{2\,c_0\,A}\,\Delta(w)\,\,\, (\mbox{ by } \eqref{equ-1})
			\end{eqnarray}
			and hence we have
			$$	l(\gamma[x,w])\le C\, \Delta(w), \,\,\,\,\, \mbox{ where } C=2e^{2\,c_0A}.$$
			Theorefore, we have
			\begin{equation}\label{equ-4}
				\min\{l(\gamma[x,w]), l(\gamma[w,y])\}\le C\Delta(w)\,\,\,\,\, \mbox{ for all } w\in \gamma[x_i,x_{i+1}].
			\end{equation}
			
			\item[\textbf{Case 2.}] Suppose $w\in \gamma[x_n,y_m]$. Then it is easy to see that
			$$\log\frac{\Delta(x_n)}{\Delta(w)}\le k_D(w,x_n)\le l_k(\gamma[w,x_n])\le l_k(\gamma[x_n,y_m]) \le c_0\,k_D(x_n,y_m) \le c_0\,A,$$
			which implies
			$$\Delta(x_n)\le e^{c_0\,A}\, \Delta(w).$$
			An easy computation shows that
			
\begin{eqnarray*}
l(\gamma[x,w])&\le& \sum_{i=0}^{n-1}l(\gamma[x_i,x_{i+1}])+l(\gamma[x_n,y_m])\\ \nonumber &\le& \sum_{i=0}^{m-1}e^{c_0\,A}\Delta(x_i)+e^{c_0\,A}\Delta(x_n)\\ \nonumber &\le& 3\,e^{c_0\, A}\Delta(x_n)\\ \nonumber
&\le& 3\,e^{2\,c_0\,A}\,\Delta(w).
\end{eqnarray*}
			
That is,
$$l(\gamma[x,w])\le C\, \Delta(w), \,\,\,\,\, \mbox{ where } C=3\,e^{2\,c_0\,A}$$
and hence we have

\begin{equation}\label{equ-5}
				\min\{l(\gamma[x,w]), l(\gamma[w,y])\}\le C\Delta(w)\,\,\,\,\, \mbox{ for all } w\in \gamma[x_n,y_m].
			\end{equation}
			\item[\textbf{Case 3.}] Suppose $w\in \gamma[y_{j+1},y_{j}]$ for some $0\le j\le m-1$. Using the same argument as in earlier case, we obtain
			$$l(\gamma[w,y])\le C\, \Delta(w), \,\,\,\,\, \mbox{ where } C=2\,e^{2c_0A}$$
			and hence
			\begin{equation}\label{equ-6}
				\min\{l(\gamma[x,w], l(\gamma[w,y]))\}\le C\Delta(w)\,\,\,\,\, \mbox{ for all } w\in \gamma[y_j,y_{j+1}].
			\end{equation}
		\end{itemize}
		Thus, combining \eqref{equ-4}, \eqref{equ-5} and \eqref{equ-6}, we conclude that $\gamma$ is a $C$- cone arc, where $C=3\,e^{2c_0A}$.
	\end{pf}
	
	\begin{thm}\label{thm-2}
		Let $D$ be a $c$- John domain in a Banach space $E$ such that $(D,k_D)$ is Gromov hyperbolic space. Fix $x,y\in D$ and let $\gamma$ be a $c_0$ neargeodesic joining $x,y\in D$. Then there exists a constant $A\ge 1$ depending only $c$, $c_0$ and $\delta$ such that the following holds:  If there are $u,v \in \gamma$ such that each $w\in \gamma{[u,v]}$ satisfies $\Delta(w)\le 2\min\{\Delta(u), \Delta(v)\}$, then $$k_D(u,v)\le A.$$
	\end{thm}
	
	\begin{pf}
		Let $D$ be a $c$-John domain in a Banach space $E$ such that $(D,k_D)$ is Gromov $\delta$-hyperbolic. By definition $(D,k_D)$ is an intrinsic space. Fix $x,y\in D$ and $c_0>1$. In view of Theorem D, let $\gamma$ be a $c_0$-neargeodesic joining $x$ and $y$. In view of Remark \ref{rem-1}, it is clear that a $c_0$-neargeodesic $\gamma$ is $(c_0,0)$- quasi-isometry. Since $D$ is a $c$-John, there is a $c$-cone arc $\alpha$ joining $x$ and $y$. We want to use Theorem F, for this we will prove the following lemma
		
		\begin{lem}\label{lem}
		Let $D$ be a $c$-John domain in a Banach space $E$ and let $\alpha$ be a $c$-cone arc joining $x,y\in D$. Let $x_0$ be a midpoint of $\alpha$, then for any $z_1,z_2 \in \alpha[x,x_0]$, we have 
	$$l_k(\alpha[z_1,z_2])\le3c\,k_D(z_1,z_2)+3c\,\log3c,$$
		and hence the quasihyperbolic arc length parametrization of $\alpha[x,x_0]$ is $(\lambda,\mu)$-quasi-hyperbolic quasi-isometry with $\lambda=3c$ and $\mu=3c\,\log(3c)$. 
		\end{lem}
		
		\begin{pf}[Proof of Lemma 3.9]
		Let $z_1,z_2\in \alpha[x,x_0]$. Without loss of generality we can assume that $z_2\in \alpha[z_1,x_0]$. Since $\alpha$ is a $c$-cone arc and $x_0$ is a mid point, for any $w\in \alpha[z_1,z_2]$, we have
		\begin{equation}\label{eq-1}
	l(\alpha[z_1,w])\le l(\alpha[x,w])\le \min \{l(\alpha[x,w]),l(\alpha[w,y])\}\le c \Delta(w),
		\end{equation}

			Next, our aim is to show that $\Delta(z_1)\le 2c\Delta(w)$. To prove this, we only need to consider the case when $d(w,z_1)\ge \Delta(z_1)/2$ because the case $d(w,z_1)<\Delta(z_1)/2$ is trivial. If $d(w,z_1)\ge \Delta(z_1)/2$, then, since $\alpha$ is a $c$-cone arc, we have
			$$c\Delta(w)\ge l(\alpha[z_1,w])\ge d(z_1,w)\ge \Delta(z_1)/2$$
			and hence 
			\begin{equation}\label{eq-2}
				\Delta(z_1)\le 2c\Delta(w)
			\end{equation}
	
	By adding \eqref{eq-1} and \eqref{eq-2}, we obtain
$$l(\alpha[z_1,w])+\Delta(z_1)\le 3c \Delta(w)$$
which gives,

\begin{equation}\label{eq-3}
\frac{1}{\Delta(w)}\le \frac{3c}{	l(\alpha[w,z_1])+\Delta(z_1)}
\end{equation}
A simple computation shows that
		
			\begin{eqnarray*}\label{3.13}
			l_k(\alpha[z_1,z_2])&\le&\int_{\alpha[z_1,z_2]}\frac{ds}{\Delta(w)}\\[2 mm] \nonumber
			&\le&\int_{0}^{l(\alpha[z_1,z_2])}\frac{3c}{t+\Delta(z_1)}\,dt \,\,\,\, (\mbox{ by } \eqref{eq-3}) \\[2 mm]\nonumber	
			&\le& 3c(\log(l(\alpha[z_1,z_2])+\Delta(z_1))-\log(\Delta(z_1))\\[2 mm] \nonumber
				&\le& 3c\,\log\left(1+\frac{l(\alpha[z_1,z_2])}{\Delta(z_1)}\right)\\[2 mm] \label{3.14}
				&\le& 3c\,\log\left(1+\frac{c\Delta(z_2)}{\Delta(z_1)}\right)\label{6.15}\\[2 mm] \nonumber
				&\le& 3c\left(\log\left(\frac{\Delta(z_2)}{\Delta(z_1)}\right)+\log\,3c\right).\\[2 mm]\nonumber
				l_k(\alpha[z_1,z_2])&\le&3c\,k_D(z_1,z_2)+3c\,\log\,3c
			\end{eqnarray*}
			Theorefore, it is clear that the quasihyperbolic arc length parametrization of $\alpha[x,x_0]$  is $(3c,3c\,\log(3c))$-quasi-isometry. This completes the proof.
		\end{pf}
		
		Since $\gamma$ is $(c_0,0)$-quasi-isometry and $\alpha[x,x_0]$ is $(3c,3c\,\log(3c))$-quasi-isometry, in view of Theorem F, we have that quasihyperbolic Hausdorff distance satisfies the following inequality
		$$k_{D}^{\mathcal{H}}(\gamma,\alpha[x,x_0])\le R=R(c,c_0,\delta).$$
		Let $u,v \in \gamma$ such that each $w\in \gamma{[u,v]}$ satisfies $\Delta(w)\le 2\min\{\Delta(u), \Delta(v)\}$. Our aim is to show that there exists a constant $A\ge 1$ such that $k_D(u,v)\le A$. For this we need to estimate $k_D(u,v)$. Since $k_{D}^{\mathcal{H}}(\gamma,\alpha[x,x_0])\le R$ {\it i.e.}, 
		$$\max\left\{\sup_{x\in \gamma}k_D(x,\alpha[x,x_0]), \sup_{y\in \alpha[x,x_0]}k_D(\gamma,y)\right\}\le R$$
		Theorefore, 
		$$\sup_{x\in \gamma}k_D(x,\alpha[x,x_0])\le R$$
		In particular,
		$$M_1=k_D(u,\alpha[x,x_0])=\inf_{b\in \alpha[x,x_0]}k_D(u,b)\le R$$
		and 
		$$M_2=k_D(v,\alpha[x,x_0])=\inf_{b\in \alpha[x,x_0]}k_D(v,b)\le R$$
		There exist sequences $\{a_i\}$ and $\{b_i\}$ such that $k_D(u,a_i)$ converges to $M_1$ and $k_D(v,b_i)$ converges to $M_2$. Since $\alpha[x,x_0]$ is compact, $a_i$ and $b_i$ have convergent subsequence which converges to say $Z_u$ and $Z_v$ respectively such that
		$k_D(u,Z_u)\le R$ and $k_D(v,Z_v)\le R$.
		Therefore, the following inequality
		$$\log \frac{\Delta(u)}{\Delta(Z_u)}\le k_D(u,Z_u)\le R$$
		gives us that $e^{-R}\Delta(u)\le \Delta(Z_u)$. Similarly, we have $\Delta(Z_v)\le e^{R}\Delta(v)$.
		Since each $w\in \gamma{[u,v]}$ satisfies $\Delta(w)\le 2\min\{\Delta(u), \Delta(v)\}$, therefore $\Delta(v)\le2\Delta(u)$ and $\Delta(u)\le 2\Delta(v)$. This gives us that $\Delta(Z_v)\le 2\,e^{2R}\Delta(Z_u)$. Hence it easily follows that
	
\begin{equation}\label{3.15}
	k_D(Z_u,Z_v)\le l_k(\alpha[Z_u,Z_v])\le 3c\,\log\left(1+\frac{c\Delta(Z_v)}{\Delta(Z_u)}\right)
\le 3c\,\log\left(1+2\,c\,e^{2R}\right).
\end{equation}

		Therefore, we have
		\begin{eqnarray*}
			k_D(u,v)&\le& k_D(u,Z_u)+k_D(Z_u,Z_v)+k_D(Z_v,v)\\ \nonumber
			&\le& 2R+k_D(Z_u,Z_v)\\ \nonumber
			&\le& 2R+3c\,\log(1+2c\,e^{2R})\,\,\,\,\,\,\, (\mbox{by } \eqref{3.15}).
		\end{eqnarray*}
		This completes the proof of Theorem \ref{thm-2} by taking $A=2R+3c\,\log(1+2ce^{2R})$.
	\end{pf}
	
	\subsection*{ Proof of Theorem \ref{main}} Let $D$ be a $c$- John domain such that $(D,k_D)$ is $\delta$-hyperbolic for some $\delta$ and let $\gamma$ be a $c_0$- neargeodesic in $D$. Then from Theorem \ref{thm-2} we have that if there are $u,v \in \gamma$ such that each $w\in \gamma{[u,v]}$ satisfies $\Delta(w)\le 2\min\{\Delta(u), \Delta(v)\}$, then 
$k_D(u,v)\le A.$ Therefore, by Theorem \ref{thm-1}, we conclude that $\gamma$ is a $b$-cone arc, where $b$ depends only on $c,\delta$ and $c_0$. This completes the proof of Theorem \ref{main}.\\
	
	\noindent\textbf{Acknowledgement:} We sincerely thank the referee for the careful reading of the manuscript and identification of several mistakes. The insightful recommendations provided by the referee enabled us to significantly enhance the paper's exposition clarity. The first named author thanks SERB-CRG, and the second named author thanks PMRF-MHRD (Id: 1200297), Govt. of India for their support.


\begin{thebibliography}{99}
		\bibitem{Matti-1986}
		{\sc G. D. Anderson}, {\sc M. K. Vamanamurthy} and {\sc M. Vuorinen}, Dimension-free quasiconformal deformation in $n$-space, {\it Trans. Amer. Math. Soc.} {\bf 297} (1986), 687--706.
		
		\bibitem{Koskela-2001}
		{\sc M. Bonk}, {\sc J. Heinonen} and {\sc P. Koskela}, Uniformizing Gromov hyperbolic domains, {\it Asterisque} {\bf 270} (2001), 1--99
		
		\bibitem{Gehring-1976}
		{\sc F. W. Gehring} and {\sc B. P. Palka}, Quasiconformally Homogeneous Domains, {\it J. Anal. Math.} {\bf 30} (1976), 172--199.
		
		\bibitem{Gehring-1979}
		{\sc F. W. Gehring} and {\sc B. G. Osgood}, Uniform domains and the Quasi-Hyperbolic Metric, {\it J. Anal. Math.} {\bf 36} (1979), 50--74.
		
		\bibitem{Gehring-1989}
		{\sc F. W. Gehring}, {\sc K. Hag} and {\sc O. Martio}, Quasihyperbolic geodesics in John domains, {\it Math. Scand.} {\bf 65} (1989), 75--92.
		
		\bibitem{Heinonen-1989}
		{\sc J. Heinonen}, Quasiconformal mappings onto John domains, {\it Rev. Math. Iber.} {\bf 5} (1989), 97--123.
		
		\bibitem{John-1961}
		{\sc F. John}, Rotation and strain, {\it Pure Appl. Math.} {\bf 14} (1961), 391--413.
		
		\bibitem{Rasila-2017}
		{\sc R. Kl\'en}, {\sc A. Rasila} and {\sc J. Talponen}, On the smoothness of quasihyperbolic balls, {\it Ann. Acad. Sci. Fenn. Math.} \textbf{42} (2017), 439--452.
		
		\bibitem{Li-2014}
		{\sc Y. Li}, Neargeodesics in John domains in Banach spaces, {\it Int. J. Math.} {\bf 25}, 1450041 (2014), https://doi.org/10.1142/S0129167X14500414. 
		
		\bibitem{Martin-1985}
		{\sc G. J. Martin}, Quasiconformal and bi-Lipschtiz homeomorphisms, uniform domains and the quasihyperbolic metric, {\it Trans. Amer. Math. Soc.} {\bf 292}(1985), 169--191.
		
		\bibitem{Martio-2011}
		{\sc Olli Martio} and {\sc J. V\"ais\"al\"a}, Quasihyperbolic geodesic in convex domains II, {\it Pure Appl. Math. Q.} \textbf{7} (2011), 395--409.
		
		\bibitem{Vaisala-1990}
		{\sc J. V\"ais\"al\"a}, Free quasiconformality in Banach spaces I, {\it Ann. Acad. Sci. Fenn. Ser. A I Math.} {\bf 15} (1990), 355--379.
		
		\bibitem{Vaisala-1991}
		{\sc J. V\"ais\"al\"a}, Free quasiconformality in Banach spaces II, {\it Ann. Acad. Sci. Fenn. Ser. A I Math.} {\bf 16} (1991), 255--310.
		
		\bibitem{Vaisala-1999}
		{\sc J. V\"ais\"al\"a}, The Free quasiworld, Quasiconformal and related maps in Banach spaces, Banach Center Publ. \textbf{48} (1999), 55--118.
		
		\bibitem{Vaisala-2004}
		{\sc J. V\"ais\"al\"a}, Broken tube in Hilbert spaces {\it Analysis} \textbf{24} (2004), 227--238.
		
		\bibitem{Vaisala-2005}
		{\sc J. V\"ais\"al\"a}, Quasihyperbolic geodesic in convex domains, {\it Results Math.} \textbf{48} (2005), 184--195.
		
		\bibitem{Vaisala-GH}
		{\sc J. V\"ais\"al\"a}, Gromov hyperbolic spaces, {\it Expo. Math.} \textbf{23} (2005), 187--231.
		
		\bibitem{Vaisala-GH-1}
		{\sc J. V\"ais\"al\"a}, Hyperbolic and Uniform domains in Banach spaces, {\it Ann. Acad. Sci. Fenn.} \textbf{30} (2005), 261--302.
		
		
		
		
		
		\bibitem{Ponnusamy-2022}
		{\sc Q. Zhou} and {\sc S. Ponnusamy}, Quasihyperbolic geodesics are cone arcs, {\it J. Geom. Anal.} {\bf 34}, 2 (2024). https://doi.org/10.1007/s12220-023-01448-x.
		
		\bibitem{Ponnusamy-2024-I}
		{\sc Q. Zhou} and  {\sc S. Ponnusamy},  Gromov hyperbolic John is quasihyperbolic John I, {\it Ann. Sc. Norm. Super. Pisa cl. Sci.} (2024), 19 pages;  DOI:10.2422/2036-2145.202207006.
		
		\bibitem{Rasila-2022}
		{\sc Q. Zhou}, {\sc Y. Li} and {\sc A. Rasila}, Gromov hyperbolicity, John Spaces, and Quasihyperbolic Geodesics, {\it J. Geom. Anal.} {\bf 32}, 228 (2022). https://doi.org/10.1007/s12220-022-00968-2.
		
		\bibitem{Matti-book}
		{\sc M. Vuorinen}, Conformal geometry and Quasiregular mappings, Lecture notes in mathematics, Vol. 1319 ({\it Springer, 1988}).
	\end{thebibliography}
\end{document}